\documentclass{amsart}
\usepackage{amssymb}
\usepackage{url}
\usepackage{algorithmic}

\newcommand\blfootnote[1]{%
  \begingroup
  \renewcommand\thefootnote{}\footnote{#1}%
  \addtocounter{footnote}{-1}%
  \endgroup
}

\newcommand{\re}{\mathrm{Re}}% partie reelle
\newcommand{\im}{\mathrm{Im}}% partie imaginaire
\newtheorem{theorem}{Theorem}[section]

\theoremstyle{definition}

\theoremstyle{remark}

\numberwithin{equation}{section}

\begin{document}

\title{On the Power Series Expansion of the Reciprocal Gamma Function}

\author[L. Fekih-Ahmed]{Lazhar Fekih-Ahmed}
\address{\'{E}cole Nationale d'Ing\'{e}nieurs de Tunis, BP 37,
Le Belv\'{e}d\`{e}re 1002 , Tunis, Tunisia}

\curraddr{} \email{lazhar.fekihahmed@enit.rnu.tn}
\thanks{}

%    \subjclass is required.


\date{July 14, 2014}

\dedicatory{}

%    Abstract is required.
\begin{abstract}
Using the reflection formula of the Gamma function, we derive a
new formula for the Taylor coefficients of the reciprocal Gamma
function. The new formula provides  effective asymptotic values
for the coefficients even for very small values of the indices.
Both the sign oscillations and the leading order of growth are
given.
\end{abstract}

\maketitle

\bibliographystyle{amsplain}

%%%%%%%%%%%%%%%%%%%%%%%%%%%%%%%%%%%%%%%%%%%%%%%%%%%%%%%%%%%%%%%%%%%%%%%%

%    Templates for common elements of a journal article; for additional
%    information, see the AMS-LaTeX instructions manual, instr-l.pdf,
%    included in the MCOM author package, and the amsthm user's guide,
%    linked from http://www.ams.org/tex/amslatex.html .

%    Section headings
\section{Introduction}\label{sec1}

The reciprocal Gamma function is an entire function with a Taylor
series given by

\begin{equation}\label{sec1-eq1}
\frac{1}{\Gamma(z)}=\sum_{n=1}^{\infty}a_{n}z^{n}.
\end{equation}

It has been a challenge since the time of Weierstrass to compute
or at least estimate  the coefficients of the reciprocal Gamma
function. The main reason is the ubiquitous presence of the
reciprocal gamma function in analytic number theory and its
various connections to other transcendental functions (for example
the Riemann zeta function). Since Bourguet
\cite{bourguet:eulerienne} who was the first to calculate the
first 23 coefficients, there has been very few publications, to
the author's knowledge, on accurate calculations of the
coefficients beyond $a_{50}$.

Knowing that an effective asymptotic formula is always useful as
an independent check for the sign and value of the coefficients
for very large values of $n$, it is important to have such a
formula in order to enhance the  calculations. The only asymptotic
formula that is known to date is that of Hayman
\cite{hayman:stirling}.

With regard to the computation of the coefficients of the
reciprocal Gamma function, there are basically three known methods
\cite{wrench:gamma, bourguet:eulerienne, bornemann:cauchy}. The
first method is due to Bourguet \cite{bourguet:eulerienne}. It
consists in exploiting the recursive formula

\begin{equation}\label{sec1-eq2}
na_{n}=\gamma a_{n-1}-\zeta(2)a_{n-2}+\zeta(3)a_{n-3}-\cdots
+(-1)^{n+1}\zeta(k); n> 2,
\end{equation}

with $a_1=1$, $a_2=\gamma$, the Euler constant, and $\zeta(k)$ is
the zeta function of Riemann.

It has been noticed in \cite{bornemann:cauchy} that this method
suffers from severe numerical instability; all digits are lost
from $n\ge 27$.

The second method is based on Cauchy's formula for the
coefficients of Taylor series using circular contours:

\begin{align}
a_n&=\frac{1}{2\pi i}\int_{|z|=r}\frac{1}{z^{n+1}\Gamma(z)} \,dz\nonumber\\
&=\frac{1}{2\pi
r^n}\int_{0}^{2\pi}\frac{e^{-in\theta}}{\Gamma(re^{i\theta})}
\,d\theta,\label{sec1-eq3}
\end{align}

where $r$ can be between 0 and $\infty$ since the reciprocal Gamma
function is entire.

The integral (\ref{sec1-eq3}) can be evaluated with the many
existing quadrature rules such as the trapezoidal rule or the
Gauss-Legendre quadrature. Particular attention is given to the
method discovered by Lyness \cite{lyness:quadrature} which uses
the trapezoidal rule in conjunction with the discrete Fourier
transform. It is very fast and provides good results
\cite{trefethen:trapezoidal} as long as the radius of the contour
is properly selected.

Although the radius $r$ of the contour can  theoretically be
arbitrarily chosen, the effects of the value of $r$ on
approximation and round-off errors are numerically very different.
A comprehensive investigation for choosing a good radius $r$ has
been carried out in \cite{bornemann:cauchy}, where it has been
shown that as $n$ increases so does the good $r$.

In \cite{schmelzer:trefethen}, different quadrature formulas,
using also the method of contour integration, have been
investigated for the calculation of $a_n$. The contour chosen is
no longer circular but chosen as the Hankel contour. The
reciprocal Gamma function is represented using Heine's formula
\cite{heine:gamma}:

\begin{equation}\label{sec1-eq4}
\frac{1}{\Gamma(z)}= \frac{1}{2\pi
i}\int_{\mathcal{C}}e^{t}t^{-z}\,dt,
\end{equation}

where $\mathcal{C}$ consists of the three parts $C=C_{+}\cup
C_{\epsilon}\cup C_{-}$:  a path which extends from
$(\infty,\epsilon)$, around the origin counter clockwise on a
circle of center the origin and of radius $\epsilon$ and back to
$(\epsilon,\infty)$, where $\epsilon$ is an arbitrarily small
positive number.

Lastly, the third  method for calculating the coefficients of the
reciprocal Gamma function for large values of $n$ is not a
numerical one. It consists in approximating the coefficients using
an  asymptotic formula. The first attempt was initiated by
Bourguet \cite{bourguet:eulerienne} who found the following upper
bound

\begin{align}
a_n&\le \frac{(-1)^n}{\pi \Gamma(n+1)}\frac{e \pi^{n+1}}{n+1}+\frac{4}{\pi^2\sqrt{\Gamma(n+1)}}\nonumber\\
&\lesssim  \frac{4}{\pi^2\sqrt{\Gamma(n+1)}}.\label{sec1-eq5}
\end{align}

But the first  systematic study to obtain an asymptotic formula
for the coefficients was carried out by Hayman
\cite{hayman:stirling} theoretically, and by Bornemann
\cite{bornemann:cauchy} numerically (see also
\cite{berry:oscillations} for the related phenomenon of
oscillations of the derivatives).

In this paper, we will give  a new effective asymptotic formula
for the coefficients $a_n$. With the formula, we obtain the sign
oscillations and the leading order of growth of the coefficients.
We will show that our results can be considered very accurate even
for very small values of $n$.

\section{An Integral Formula For The Coefficients $a_n$}\label{sec2}

Let's replace $z$ by $z-1$ into the series (\ref{sec1-eq1}), we
have

\begin{equation}\label{sec2-eq1}
\frac{1}{\Gamma(z-1)}=a_1(z-1)+a_2(z-1)^2+a_3(z-1)^3\cdots ,
\end{equation}

and dividing both sides by $z-1$, we get

\begin{equation}\label{sec2-eq2}
\frac{1}{\Gamma(z)}=\frac{1}{z-1}\left[
a_1(z-1)+a_2(z-1)^2+a_3(z-1)^3 \cdots \right]
\end{equation}

To obtain an integral formula for the reciprocal Gamma function,
we start from Euler's reflection formula

\begin{equation}\label{sec2-eq3}
\Gamma(z)\Gamma(1-z)=\frac{\pi}{\sin(\pi z)}
\end{equation}

to get

\begin{equation}\label{sec2-eq4}
\frac{1}{\Gamma(z)}=\frac{\sin(\pi z)}{\pi}\Gamma(1-z).
\end{equation}

Now, for $\re(z)<2$, we can write
\begin{align}
\frac{1}{\Gamma(z)}&=\frac{\sin(\pi z)}{\pi}\Gamma(1-z)\nonumber\\
&=\frac{\sin(\pi (z-1))}{\pi(z-1)}\Gamma(2-z)\nonumber\\
&=\frac{\sin(\pi (z-1))}{\pi(z-1)}\int_{0}^{\infty}e^{-t}
t^{1-z}\,dt.\label{sec2-eq5}
\end{align}

By observing that  $\sin(\pi(z-1))=\frac{e^{i\pi (z-1)}-e^{-i\pi
(z-1)}}{2i}$, we can rewrite (\ref{sec2-eq5}) as

\begin{equation}\label{sec2-eq6}
\frac{1}{\Gamma(z)}=\frac{1}{z-1}\frac{1}{2\pi i}
\int_{0}^{\infty}e^{-t}
\big[e^{(z-1)(-\log(t)-i\pi)}-e^{(z-1)(-\log(t)+i\pi)}\big]\,dt.
\end{equation}

And if we compare the two equations (\ref{sec2-eq6}) and
(\ref{sec2-eq2}), we deduce that the coefficients $a_n$ for $n\ge
1$ are given by

\begin{eqnarray}\nonumber
a_n&=& \frac{1}{2\pi i n!}\int_{0}^{\infty}e^{-t}\lim_{z\to
1}\frac{d^n}{d z^n} \left \{ e^{(z-1)(-\log(t)+i\pi)}-e^{(z-1)(-\log(t)-i\pi)}\right \}\,dt\nonumber\\
&=&\frac{1}{2\pi i n!}\int_{0}^{\infty}e^{-t} \left \{
(-\log(t)+i\pi)^n-(-\log(t)-i\pi)^n\right
\}\,dt\nonumber\\
&=&\frac{1}{\pi n!}\int_{0}^{\infty}e^{-t} \im\left \{
(-\log(t)+i\pi)^n\right \}\,dt,\label{sec2-eq7}
\end{eqnarray}

where $\im$ stands for the imaginary part. This is our expression
of the coefficients $a_n$, described in the following

\begin{theorem}\label{sec2-thm2}
The coefficients $a_{n}$ are given by
\begin{equation}\label{sec2-eq8}
a_n=\frac{(-1)^n}{\pi n!}\int_{0}^{\infty}e^{-t}\im
\big\{(\log{t}-i\pi)^{n}\big\}\,dt.
\end{equation}
\end{theorem}

Theorem~\ref{sec2-thm2}  permits an effective asymptotic evaluation of
the constants $a_{n}$\footnote{The theorem is almost evident and
easy to derive. It is hard to believe  that it has not been
discovered before. To the author's knowledge, the integral formula
(\ref{sec2-eq8}) is new and seems to be nonexistent in the
literature.}. It is the subject of the next section.

\section{Asymptotic Estimates of the Coefficients}\label{sec3}

This section is dedicated to approximating the complex-valued
integral

\begin{equation}\label{sec3-eq1}
I(n)=\int_{0}^{\infty}e^{-t}(\log{t}-i\pi)^{n}\,dt
\end{equation}

using the saddle-point method \cite{copson:1965,debruijn:1981}. By
the change of variables $t=nz$, our integral becomes

\begin{align}
I(n)&=n\int_{0}^{\infty}
e^{-nz}\left\{\log\left(nz\right)-i\pi\right\}^{n}\,dz\nonumber\\
&=n\int_{0}^{\infty} e^{n\left\{-z+
\log\left[\log\left(nz\right)-i\pi\right]\right\}}\,dz.\label{sec3-eq2}
\end{align}

If we define

\begin{equation}\label{sec3-eq3}
f(z)=-z+ \log\left(\log\left(nz\right)-i\pi\right),
\end{equation}

then the saddle-point method consists in  deforming the path of
integration into a path which goes through a saddle-point at which
the derivative $f^{\prime}(z)$, vanishes. If $z_0$ is the
saddle-point at which the real part of $f(z)$ takes the greatest
value, the neighborhood of $z_0$ provides the dominant part of the
integral as $n\to \infty$ \cite[p. 91-93]{copson:1965}. This
dominant part provides an approximation of the integral and it is
given by the formula

\begin{align}
I(n)&\approx
ne^{nf(z_0)}\left(\frac{-2\pi}{nf^{\prime\prime}(z_0)}
\right)^{\frac{1}{2}}.\label{sec3-eq4}
\end{align}

In our case, we have

\begin{align}
f^{\prime}(z)&=-1+\frac{1}{z\left(\log\left(nz\right)-i\pi\right)},\quad\text{and}
\label{sec3-eq5}\\
f^{\prime\prime}(z)&=\frac{-1}{z^2\left(\log\left(nz\right)-i\pi\right)}-
\frac{1}{z^2\left(\log\left(nz\right)-i\pi\right)^2}.\label{sec3-eq6}
\end{align}

The saddle-point $z_0$ should verify the equation

\begin{align}
&\quad z_0\left(\log\left(nz_0\right)-i\pi\right)=1\nonumber\\
\Leftrightarrow&\quad nz_0e^{-i\pi}\log\left(nz_0
e^{-i\pi}\right)= n e^{-i\pi}.\label{sec3-eq7}
\end{align}

The last equation is of the form $v\log{v}=b$ whose solution can
be explicitly written using the branch $k=-1$ of the Lambert
$W$-function\footnote{The principal branch of the Lambert
$W$-function is denoted by $W_0(z)=W(z)$. The principal branch
$W_0(z)$ and the branch $W_{-1}(z)$ are the only branches of $W$
that take on real values. The other branches of $W$ have the
negative real axis as the only branch cut closed on the top for
counter clockwise continuity. In our equation (\ref{sec3-eq7}),
the argument is $-\pi$ and not $\pi$ and so the solution belongs
to the branch  of $W_{-1}$. See \cite{corless:lambert} for a
excellent discussion and explanation of all the branches of $W$.}
\cite{corless:lambert}:

\begin{align}
v=e^{W_{-1}(b)}.\label{sec3-eq8}
\end{align}

The saddle-point solution to our equation (\ref{sec3-eq7}) is
given by

\begin{align}
z_0=\frac{e^{-i\pi}}{n}e^{W_{-1}\left(ne^{-i\pi}\right)}=\frac{e^{W_{-1}}(-n)}{-n},\label{sec3-eq9}
\end{align}

and at the saddle-point, we have the values

\begin{align}\label{sec3-eq10}
f(z_0)&=-z_0-\log{z_0}\\
f^{\prime\prime}(z_0)&=-1-\frac{1}{z_0}.\label{sec3-eq11}
\end{align}

Therefore, the saddle-point approximation of our integral
(\ref{sec3-eq1}) is given by

\begin{equation}\label{sec3-eq12}
I(n)\approx\sqrt{2\pi n}e^{-nz_0}\frac{z_0^{\frac{1}{2}-n}}{
\sqrt{1+z_0}}.
\end{equation}

Now since $a_n= \frac{(-1)^n}{\pi n!}\im\left\{I(n)\right\}$, we
arrive at our main result:

\begin{theorem}\label{sec3-thm1}
Let $z_0=\frac{e^{W_{-1}}(-n)}{-n}$, where $W_{-1}$ is the branch
$k=-1$ of the Lambert $W$-function. For $n$ large enough, the
Taylor coefficients of the reciprocal Gamma function can be
approximated by
\begin{equation}\label{sec3-eq13}
a_n\approx  (-1)^n \sqrt{\frac{2}{\pi}}\frac{\sqrt{
n}}{n!}\im\left\{e^{-nz_0}\frac{z_0^{\frac{1}{2}-n}}{
\sqrt{1+z_0}}\right\}.
\end{equation}
\end{theorem}

Bornemann's derivation \cite{bornemann:cauchy} of Hayman's
asymptotic formula for the coefficients $a_n$ is given
by\footnote{The formula of Bornemann differs from that of Hayman
in the phase approximation. The orignal approximation given by
Hayman is
$\phi_n=\left(n-\frac{1}{2}\right)\left(\frac{\sin^2\theta_n}{\theta_n}-\theta_n\right)$.
For the calculations, both phase approximations give essentially
the same results.}

\begin{align}\label{sec3-eq14}
a_n&\sim\frac{\sqrt{2}}{\pi n}\frac{1}{|\Gamma\left(r_{n}
e^{i\theta_n}\right)|r_n^{n}}\cos\phi_n,
\end{align}
where

\begin{align}\label{sec3-eq15}
z_n&=r_ne^{i\theta_n}=e^{W\left(\frac{1}{2}-n\right)}\\
\phi_n&=\left(n-\frac{1}{2}\right)\left(\frac{\sin^2\theta_n}{\theta_n}-\theta_n\right)-
\frac{1}{2}\left(\cot\theta_n-
\theta_n\csc^2\theta_n\right).\label{sec3-eq16}
\end{align}

Note that both formulas use the lambert $W$-function. Our formula
will be compared to Hayman's formula in the next section.

We can also find an asymptotic formula of our $a_n$ as a function
of $n$ only by resorting to the following asymptotic development
of the  branch of $W_{-1}(z)$ \cite{corless:lambert}:

\begin{equation}\label{sec3-eq17}
W_{-1}(z)=\log(z-2\pi i)-\log\left(\log{z-2\pi i}\right)+\cdots
\end{equation}

For $n\gg 1$ we can write

\begin{align}\label{sec3-eq18}
z_0&\sim \frac{-n-2\pi i}{-n\log(-n-2\pi i)}\sim
\frac{1}{\log{n}-\pi i} \sim \frac{e^{-i
\arctan\left(\frac{\pi}{\log{n}}\right)}}{\sqrt{(\log{n})^2+\pi^2}
}\sim \frac{e^{i\frac{\pi}{\log{n}}}}{ \log{n} },
\end{align}

\begin{align}\label{sec3-eq19}
\frac{z_0^{\frac{1}{2}-n}}{ \sqrt{1+z_0}}\sim
\frac{1}{z_0^{n}}&\sim \frac{(\log{n})^{n}}{e^{i
\frac{n\pi}{\log{n}}}},
\end{align}

and

\begin{align}\label{sec3-eq20} e^{-nz_0}&\sim
e^{\frac{-n}{\log{n}}e^{-i\frac{\pi}{\log{n}}}}\sim
e^{\frac{-n}{\log{n}}} ,
\end{align}

and using Stirling formula $n!\sim \sqrt{2\pi n}\big
(\frac{n}{e}\big )^n$, this  yields the second approximation

\begin{align}\label{sec3-eq21}
a_n&\sim \frac{(-1)^{n+1}}{\pi}
e^{-n\log{n}+n\log{\log{n}}+n-\frac{n}{\log{n}}}\sin\left(
\frac{n\pi}{\log{n}}\right).
\end{align}

Equation (\ref{sec3-eq21}) is a rough approximation. It will not
be used for calculations. It only provides the leading order of
growth and the sign oscillations of the coefficients. However, one can  use 
the approximation (\ref{sec3-eq21}) to easily prove that the order of the reciprocal gamma function is 1 and that its type is maximal. 

\section{Numerical Results and Conclusion}\label{sec4}

We implemented the  formula of Theorem~\ref{sec3-thm1} and
Hayman's formula (\ref{sec3-eq14}) in
Maple\textsuperscript{\scriptsize{\texttrademark}}
\blfootnote{\textsuperscript{\scriptsize{\texttrademark}}Maple is
a trademark of Waterloo Maple Inc.}. For a given value of $n$, the
following code computes the value of $a_n$:

\begin{algorithmic}

\STATE \texttt{w0 := LambertW(-1, n*exp(-I*Pi)):}

\STATE \texttt{z0 := exp(-I*Pi)*exp(w0)/n:}

\STATE \texttt{f := -z0-ln(z0):}

\STATE \texttt{fpp := -1-1/z0:}

\STATE
\texttt{evalf((-1)**n*Im(n*sqrt((-2*Pi)*(1/(n*fpp)))\\
*exp(n*f))/(Pi*factorial(n)));}

\end{algorithmic}

The approximations (\ref{sec3-eq13}) and  (\ref{sec3-eq14}) were
examined and compared to the exact values for $n$ from $1$ to $20$
given in \cite{wrench:gamma}. They are displayed in
Table~\ref{sec4-table1}. Table~\ref{sec4-table2}, displays the
approximate value of $a_{n}$ and  exact values for higher values
of $n$.

\begin{table}[!h]
\centerline{
\begin{tabular}{cccc}
  \hline
  $n$ & $a_n$ &  Formula of Theorem~\ref{sec3-thm1} & Hayman's formula\\
  \hline
2& 0.577215664&  0.471315586&0.318527853\\
3&-0.655878071 &-0.634156618&-0.745580393\\
\bf 4&\bf -0.042002635 &\bf -0.024878383& \bf 0.035835755\\
5& 0.166538611 & 0.1586548367& 0.170422513\\
6&-0.042197734& -0.0422409922&-0.055165293\\
7&-0.009621971& -0.0088055266&-0.006842089\\
8& 0.007218943&  0.0070070400& 0.007791124\\
9&-0.001165167& -0.0011689459&-0.001538105\\
10&-0.000215241&-0.0002013214&-0.000162310\\
11& 0.000128050 &0.0001248855& 0.000137477\\
12&-0.000020134&-0.0000200451&-0.000025104\\
\bf 13&\bf -0.00000125&\bf -0.000001139&\bf -0.000000054\\
14& 0.000001133& 0.0000011053& 0.000001178\\
\bf 15&\bf -2.0563384.$\bf 10^{-7}$ &\bf -2.034656492.$\bf 10^{-7}$&$\bf -2.410634519.10^{-7}$\\
\bf 16&$\bf 6.11609510.10^{-9}$& $\bf  6.506886194.10^{-9}$&  $ \bf 1.201994777.10^{-8}$\\
17&5.00200764.$10^{-9}$ &4.864046460.$10^{-9}$&4.859838872.$10^{-9}$\\
18&-1.18127457.$10^{-9}$&-1.164373917.$10^{-9}$&$-1.3136121.10^{-9}$\\
19&1.043426711.$10^{-10}$ &1.043634325.$10^{-10}$&$1.3322234.10^{-10}$\\
20&7.782263439.$10^{-12}$ &7.415156531.$10^{-12}$&$5.436583518.10^{-12}$\\
\hline
\end{tabular}
}\vspace{0.5cm}
\caption{First 20  coefficients and their approximate values given by Theorem~\ref{sec3-thm1} and
formula (\ref{sec3-eq14}).}%
\label{sec4-table1} %
\end{table}

For $n=4$ Hayman's formula gives the wrong sign. It provides a
value of the coefficient with an error of $18.5\%$ for $n=15$, and
for $n=13,16$, the error is almost $96\%$. Moreover, we can see
that for at least the values of $2\le n\le 20$ Hayman's formula is
not as good an approximation to the exact value as the formula of
Theorem~\ref{sec3-thm1}.

\begin{table}[!h]
\centerline{
\begin{tabular}{cccc}
  \hline
  $n$ &$a_n$&Formula of Theorem~\ref{sec3-thm1} & Hayman's formula  \\
  \hline
30&$1.7144063219.10^{-20}$&$1.708720889.10^{-20}$& $2.072558647.10^{-20}$ \\
40&$-1.1245843492.10^{-30}$&$-1.110270738.10^{-30}$& $-1.143814145.10^{-30}$\\
50&$-1.0562331785.10^{-41}$&$-1.051407032.10^{-41}$&$-1.211991030.10^{-41}$\\
100&$6.6158100911.10^{-106}$& $6.599969140.10^{-106}$&$7.56758012.10^{-106}$\\
150&$1.1936904502.10^{-179}$&$1.193587226.10^{-179}$&$1.4412588.10^{-179}$\\
250&$-2.4488582032.10^{-343}$& $-2.446740476.10^{-343}$&$-2.8028909.10^{-343}$\\
300&$2.90203183445.10^{-431}$& $2.900143434.10^{-431}$&$3.3306712.10^{-431}$\\
800&$-2.46251758839.10^{-1431}$&$-2.460396773.10^{-1431}$ &$-2.5852781.^{-1431}$\\
1400&$-6.07622638292.10^{-2792}$& $-6.074000773.10^{-2792}$&$-6.5759375.10^{-2792}$\\
\hline
\end{tabular}
}\vspace{0.5cm} \caption{The coefficients $a_n$ for different
higher values of $n$ and their approximations using the asymptotic
formula of Hayman  (\ref{sec3-eq14}) and the asymptotic formula of
Theorem~\ref{sec3-thm1}.}%
\label{sec4-table2} %
\end{table}

From the trend of the values in Table~\ref{sec4-table1} and
Table~\ref{sec4-table2}, we  conclude that for small values of $n$
our formula outperforms Hayman's formula and that for larger
values of $n$ both formulas give the same sign  but differ
slightly in magnitude. The asymptotic formula of this paper has
the  advantage that it does not depend on the radius $r_n$ of the
circular contour, a real advancement in  estimating the
coefficients of the Taylor series of  the reciprocal Gamma
function.

As a final remark, the asymptotic formulas of this paper can of
course be used to find an  asymptotic formula for the related
constants $b_n$ defined by the power series

\begin{equation}\label{sec4-eq1}
\frac{1}{\Gamma(z)}=z(1+z)\left[b_0+b_1z+b_2z^2\cdots \right],
\end{equation}

where the coefficients $a_n$ and $b_n$ are connected by the
relation

\begin{equation}\label{sec4-eq2}
a_n=b_{n-1}+b_{n-2}; n\ge 2.
\end{equation}

\section*{Acknowledgement}
I am grateful to Dr. Fredrik Johansson for providing the exact
values in Table~\ref{sec4-table2} to 20D for $a_n$, $n\ge 30$. The
algorithm he used for the computations is described  in his thesis
\cite{johansson:thesis}.

%%%%%%%%%%%%%%%%%%%%%%%%%%%%%%

\end{document}